\documentclass[12pt]{article}
\usepackage{amsmath}
\usepackage{amssymb}
\usepackage{url}

\usepackage[T1]{fontenc}
\usepackage{newcent}
\usepackage{diagrams}
\diagramstyle{midshaft,PostScript=dvips,nohug}

\newtheorem{lemma}{Lemma}[section]
\newtheorem{prop}[lemma]{Proposition}
\newtheorem{thm}[lemma]{Theorem}
\newtheorem{cor}[lemma]{Corollary}

\newcommand{\M}{{\mathfrak M}}

\newcommand{\pps}{parallelopipeda}
\newcommand{\pgs} {parallelograms}
\newcommand{\pg}{parallelogram}

\title{Integration of 1-forms and connections}
\author{Anders Kock}
\date{  }

\begin{document}
\maketitle

\section*{Introduction}

We shall present a geometric/combinatorial  version of the following 
general wishes:  1) closed 
group-valued 1-forms locally have primitives; 2) flat (curvature-free) connections in groupoids 
locally have trivializations; 3) spaces with flat (= curvature free) 
and symmetric (= torsion free) affine connections are locally 
affine spaces. In the presentation here,  each of two last stages 
presupposes the preceding one.

For the case of affine connection in combinatorial terms (``formation 
of infinitesimal \pgs''), we solve a problem left open in \cite{SGM} p.\ 48: 
when can \pg \; formation be extended to formation of infinitesimal 
\pps? - like in an affine space? The classical answer is: when the affine 
connection is flat and symmetric. In our geometric/combinatorial 
version, it is a consequence of Theorem 
\ref{37x} below.

The solutions we give do not depend on the real numbers; they may be 
reformulated (when coordinatized) into statements about existence of  
suitable  ``formal power 
series'', without any discussion of convergence. Such reformulations  work, when 
coordinatized,  over any 
field (or even local ring) of characteristic 0. But largely, our 
exposition is coordinate free.

We shall (at least in the present  version) freely use notation 
and concepts from \cite{SDG} and \cite{SGM}. 

 Via well adapted models of synthetic differential geometry, as 
 constructed by Dubuc, (see 
\cite{SDG}), the results can be interpreted in the category of smooth manifolds in the 
classical sense (see e.g.\ \cite{LGVI}). But some of them apply in other categories, e.g.\ 
in some categories coming  from algebraic geometry. 
We shall consider the category of formal manifolds, in the sense of 
\cite{SDG} I.17. The main thing is that the objects 
$M$ which we consider come equipped with a reflexive symmetric 
relation $\sim$, (preserved by the morphisms). For schemes $M$ in 
algebraic geometry, such $\sim$ was introduced by French algebraic 
geometry (notably Grothendieck) in the 1960s, via what was called 
{\em the first neighbourhood of the diagonal, $M_{(1)}\subseteq 
M\times M$}.

Part of the notions and proofs we develop in the present paper are phrased entirely 
in  terms of this relation $\sim$ and is purely 
combinatorial.\footnote{It is worth investigating what the present 
theory of affine connections has to do with the theory of ``edge 
symmetric double groupoids with connections'' of \cite{B}.} 
But to be specific, we consider (formal) manifolds only.

We call $\sim$ the {\em (first order) neighbour} 
relation, so $x\sim y$ is read ``$x$ and $y$ are neighbours'', or 
even (first order) {\em infinitesimal} neighbours. The set of neigbours of $x$, we 
denote $\M (x)$, the (first order) {\em monad} of $x$.

 Note that the relation $\sim$ is not assumed to be 
transitive. The transitive closure of $\sim$ is an equivalence relation, $\sim 
_{\infty}$, namely $x$ and $z$ satisfy $x\sim _{\infty}z$, if 
for some 
$k\in {\mathbb N}$, we have $x\sim _{k}z$; this in turn means  that there is a 
 chain (``$k$-path'') $x\sim y_{1}\sim y_{2}\sim \ldots y_{k-1}\sim z$. 
Therefore, we call the equivalence classe {\em (infinitesimal) path 
components}. Our theory deals with such path components, or, 
equivalently, with a path connected $M$ (= equivalence class for 
$\sim _{\infty}$). The equivalence class of $x$ is denoted 
$\M_{\infty}(x)$ (the {\em $\infty$-monad} around $x$). The notion of 
`local' is, for simplicity, taken to refer to formally open subsets, 
i.e.\ subsets which are closed under the relation $\sim _{\infty}$.

\section{Group valued 1-forms}

\subsection{Basic theory of group valued 1-forms}

The following Subsection depends on the axiomatics of synthetic differential geometry; the reader 
who wants to go straight to the combinatorics, may skip this, and 
take the conclusion Proposition \ref{QLx}, and in more general form, 
Proposition \ref{quadr+x}, as an axiom.

Let $M$ be a manifold and $G$ a group 
(not necessarily commutative, multiplication denoted $*$, unit by $1$). 
Recall (from \cite{SGM}, say) that a $G$-valued 1-form is a map $\omega : M_{(1)}\to G$ with
$\omega(x,x)=1$ for all $x\in M$ (and   with 
$\omega(y,x)=\omega (x,y)^{-1}$; this can often be {\em deduced}, see 
\cite{SGM} Proposition 6.1.3). It is {\em closed} if
\begin{equation}\label{closedx}
\omega (x,y)* \omega(y,z) = \omega(x,z), 
\end{equation}
whenever $x,y$ and $z$ are mutual neighbours.
In particular, for a closed 1-form $\omega$, we have for mutual 
neighbours $x,y,z$ that  $\omega (x,y)*\omega 
(y,z)$ is independent of $y$. We may ask whether this independence of 
$y$ also applies if we do not assume that $x\sim z$. We shall prove 
in the context of synthetic differential geometry,  
that for closed forms, this independence indeed obtains, under the auxiliary assumption 
that $G$ is (isomorphic to)  a matrix group in the following sense: there exists an 
associative unitary algebra $(W,*)$, such that $G$ is a subgroup of the 
multiplicative monoid, and such that $W$ is a KL vector space, in the 
sense of \cite{SGM} 1.3. (Think 
of $(W,*)$ as a matrix algebra.) So for $x\sim y$, we have $\omega 
(x,y)\sim 1$, so it is of the form $1+d$ for some $d\in D(W)$ (= the 
set of $\sim$-neighbours of $0\in W$, or $\M(0)$).
 Therefore, 
choosing a coordinate chart $U\to M$ around $x$ and $y$, where 
$U$ is a formally open subset of a KL vector space $V$, and 
identifying  points in the image of the charts by their coordinates in the KL vector space 
$V$, the function $\omega$ may in $U$ be 
expressed in the form 
$$\omega (x,y) = 1+\Omega (x;y-x)$$
for a unique function $\Omega :U\times V \to W$, linear in 
the second argument (using the KL property). Recall that $x\sim y$ in $U$ means that 
$y-x \in D(V)$, (the first order infinitesimal neighbourhood of $0\in 
V$). Thus the relation between $\omega$ and $\Omega$ may equally be 
expressed that for $x\in U$ and $d\in D(V)$, we have
$\omega (x,x+d)= 1+\Omega (x;d)$.

The following  calculation is basically identical to some that occurs in the 
proof of
Proposition 6.2.5 in \cite{SGM} (where the present $\Omega$ is denoted 
$l\omega$, and the present $(W,*)$ is denoted $(A, \cdot )$). But note 
that in loc.cit., it is assumed that $x,y,z$ are mutual neighbours, 
whereas we here do not assume that $x\sim z$, but only that $x\sim y \sim 
z$. Thus we have $y=x+d_{1}$ and $z=y+d_{2} = x+d_{1}+d_{2}$ with 
$d_{1}$ and $d_{2}$ in $D(V)$; but we are not assuming that 
$d_{1}+d_{2}\in D(V)$.

So the $x$, $y$, and $z$ considered are of the form $x$, $x+d_{1}$, 
and $x+d_{1}+d_{2}$, respectively, with $d_{1}$ and $d_{2}$ in $D(V)$.
We calculate for such $(d_{1},d_{2}) \in D(V)\times D(V)$ the 
expression for $\omega (x,y)*\omega (y,z)$ in  terms of 
$\Omega$:  
\begin{eqnarray*}\omega(x,y)*\omega(y,z)&=&(1+\Omega (x;d_{1}))*(1+ \Omega (x+d_{1}; d_{2}))\\
&=&1+\Omega (x;d_{1})+\Omega (x+d_{1}; d_{2})+ \Omega (x;d_{1})*\Omega 
(x+d_{1}; d_{2}).\end{eqnarray*}  By Taylor expansion,   $\Omega 
(x+d_{1};d_{2}) = \Omega (x;d_{2})+d\Omega (x;d_{1},d_{2})$; 
substituting this in the two places where $\Omega 
(x+d_{1};d_{2})$ occurs, allows us to 
continue
\begin{eqnarray*}&=1+\Omega (x;d_{1})&+\Omega (x; d_{2})+d\Omega (x; d_{1},d_{2})+\\ 
&&+\; \Omega (x;d_{1})*\Omega(x; d_{2}) + \Omega (x; d_{1})*d\Omega (x; 
d_{1}, d_{2})\end{eqnarray*}
The last term here contains $d_{1}$ in a bilinear way, so it 
vanishes. So we are left with
$$1+\Omega (x;d_{1})
+\Omega (x; d_{2})+d\Omega (x; d_{1},d_{2}) 
+\Omega (x;d_{1})*\Omega(x; d_{2}),$$
so  (using $\Omega (x;d_{1})
+\Omega (x; d_{2})= \Omega (x; d_{1}+d_{2})$), we conclude
\begin{equation}\omega(x,y)*\omega(y,z)=1+\Omega (x; d_{1}+d_{2})+d\Omega (x; d_{1},d_{2}) 
+\Omega (x;d_{1})*\Omega(x; d_{2}).\label{basicx}\end{equation}



\begin{prop}\label{QLx}[Quadrangle Law] If $\omega$ is a closed $G$-valued form, then for $x\sim 
y\sim z$, we have that  $\omega (x,y)*\omega (y,z)$ is independent of 
$y$. \end{prop}
(Note that we 
cannot shortcut the conclusion of the Proposition by saying: 
``in fact, $\omega (x,y)*\omega (y,z)$ equals $\omega (x,z)$''; 
for, $\omega (x,z)$ only makes sense if $x\sim z$.)

\medskip

\noindent{\bf Proof.}  We pick a coordinate chart $U$ as above, in particular, 
$y=x+d_{1}$ and $z= x+d_{1}+d_{2}$. In terms of these coordinates, we have derived the 
expression 
(\ref{basicx}) for $\omega (x,y)*\omega (y,z)$. If $d_{1}+d_{2} \in 
D(V)$, we have \mbox{$\omega (x,z)= 
1+\Omega (x;d_{1}+d_{2})$,} so if further $\omega$ is closed, we therefore 
have, by subtracting from (\ref{basicx}), that
\begin{equation}\label{tilde}d\Omega 
(x;d_{1},d_{2})+\Omega(x;d_{1})*\Omega(x;d_{2}) =0.\end{equation}
For fixed $x$, the function $d\Omega 
(x;v_{1},v_{2})+\Omega(x;v_{1})*\Omega(x;v_{2})$ is a bilinear 
function \mbox{$V\times V \to W$.} By the equation (\ref{tilde}), this function vanishes 
when  $d_{1}$, $d_{2}$,  and $d_{1}+d_{2}$ are in $D(V)$.
We leave to the reader to prove that if $d_{1}$ and $d_{2}$ are in 
$D(V)$, 
then $d_{1}+d_{2}\in D(V)$ iff $d_{1}-d_{2} \in D(V)$, i.e.\ iff 
$d_{1}\sim d_{2}$, (use the characterization of 
$D(V)$ in terms of symmetric bilinear $V\times V \to R$, Proposition 
1.2.12 in \cite{SGM}), or again,
iff $(d_{1},d_{2})\in \tilde{D}(2,V)$ (as defined in \cite{SGM} 1.2).
So it follows 
(Proposition 1.3.3 in \cite{SGM}) that the expression in (\ref{basicx})  only depends on 
$d_{1}+d_{2}$. For $x,y,z$, this says that $\omega 
(x,y)*\omega (x,z)$ does not depend on $y$, (in coordinates: it does 
not depend on $d_{1}$), but only on $x$ and $z$. And this assertion does not depend 
on the choice of chart. This proves the Proposition.

\medskip
(The converse is also true: if $\omega (x,y)*\omega (y,z)$ is independent of 
$y$, then $\omega$ is closed. We leave this as an exercise.)

\medskip
The reason for the 
name ``quadrangle law'' is that the conclusion may expressed by 
saying that given a $\sim$-quadrangle, meaning four points 
$x,y_{1},y_{2},z$ with $x\sim y_{1}\sim z$ and $x\sim y_{2}\sim z$,  we have (for $\omega$ 
closed) that 
$\omega(x,y_{1})*\omega (y_{1},z)=\omega(x,y_{2})*\omega (y_{2},z)$. 
This equality we shall express as an equality of two ``path 
integrals'', or ``curve integrals'' of the 1-form $\omega$ along the  
periphery of the quadrangle. 

We shall, more generally, describe  path integrals of a $G$-valued 1-forms $\omega$ 
along ``paths'' of arbitrary finite length. We 
consider the formal (infinitesimal) substitute of the notion of {\em 
path} $\underline{x}$, for which the task is to describe the 
``path integral''  $\int 
_{\underline{x}}\omega \in G$:

We define an
 {\em $n$-path} $\underline{x}$ in a manifold $M$ to be an 
$n+1$-tuple \newline \mbox{$(x_{0}, x_{1}, 
\ldots , x_{n})$} of points in $M$ with
$x_{i}\sim x_{i+1}$ for $i=0, \ldots ,n-1$. The point $x_{0}$ is the 
{\em domain} of $\underline{x}$, and  the point $x_{n}$ is the {\em 
codomain} of $\underline{x}$.
If $\omega$ is a $G$-valued 1-form on $M$, we define the ``path 
integral'' $\int _{\underline{x}}\omega$ by 
$$\int _{\underline{x}}\omega : = \omega (x_{0}, x_{1})*\omega 
(x_{1},x_{2})*\ldots *\omega (x_{n-1},x_{n}).$$
So for $n=1$, $\int_{\underline{x}}\omega = \omega (x_{0},x_{1})$.

\begin{prop}\label{quadr+x}If $\omega$ is a closed $G$-valued 1-form 
on a manifold $M$, then
$\int _{\underline{x}}\omega $ only depends on the 
domain and the codomain of the path $\underline{x}$.
\end{prop}
Note that for $n=2$, this is a restatement of the Proposition \ref{QLx}.

\medskip

\noindent{\bf Proof.} As in the proof of the Proposition \ref{QLx}, we pick 
an arbitrary chart $U$ contaning all the $x_{i}$s of the path; so 
the path (say, an $n$-path) may be presented with $x_{0}$, and a sequence 
$\underline{d}=d_{1}, d_{2}, \ldots ,d_{n}$ (with $d_{i}\in D(V)$) with $x_{i}=x_{i-1}+d_{i}$
 for $i=1, \ldots ,n$. From Proposition \ref{QLx}  follows that
 $$\int_{\underline{x}}\omega = \int_{\underline{x'}}\omega,$$
 where $\underline{x'}$ is obtained from $\underline{x}$ by swapping 
 the $i$th and $(i+1)$st of the $d_{j}$s ($i=1, \ldots ,n-1$), so as 
 to obtain a new point $x_{i}'$ (this $x_{i}'$ is something 
 that depends on the chart)  \bigskip

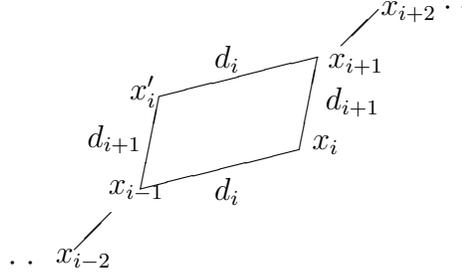
\begin{picture}(100,120)(-110,-40)
\put(20,6){\line(4,1){60}}
\put(20,6){\line(1,5){7}}
\put(80,21){\line(1,5){7}}
\put(27,41){\line(4,1){60}}

\put(8,4){$x_{i-1}$}
\put(85,21){$x_{i}$}
\put(16,40){$x_{i}'$}
\put(91,52){$x_{i+1}$}
\put(111,72){$x_{i+2}\cdot \cdot$}
\put(90,36){$d_{i+1}$}
\put(48,2){$d_{i}$}
\put(48,52){$d_{i}$}
\put(0,22){$d_{i+1}$}
\put(-5,-17){\line(1,1){14}}
\put(-30, -25){$\cdot \; \cdot $}

\put(-12,-22){$x_{i-2}$}

\put(96,60){\line(1,1){14}}



\end{picture}

\noindent We can thus swap any two conscutive entries in the 
sequence of $d_{j}$, without changing the value of the integral; and since neighbour transpositions generate the 
whole symmetric group $S_{n}$ of permutations  $\sigma$ of $n$ letters, it 
follows that (for closed $\omega$)
\begin{equation}\label{Sx}\int_{\underline{x}}\omega 
=\int_{\sigma(\underline{x})}\omega ,\end{equation}
where  $\sigma(\underline{x})$ replaces the 
$x_{i}=x_{0}+\sum_{j=1}^{i}d_{j}$ in the original $\underline{x}$ by 
$x_{i}':=x_{0}+\sum_{j=1}^{i}d_{\sigma (j)}$. 
For fixed $x_{0}$, we therefore have a map 
which is 
invariant under the $n!$ permutations of the $n$ input entries
$(d_{1}, d_{2}, \ldots ,d_{n})$, 
By the ``Symmetric Functions Property''   in its geometric 
manifestation, \cite{DK} Theorem 2.1, it follows that (\ref{Sx}), as a 
function of the $d_{i}$s, factors (in fact uniquely) across 
the addition map $D(V)^{n}\to D_{n}(V)$, i.e.\ it depends only of the 
sum $\sum d_{j}$, not on the indivual $d_{j}$'s. Equivalently, 
$\int_{\underline{x}}\omega$ only depends on $x_{0}$ and $x_{n}$. 
This is now a statement which does not mention any particular chart. 
This proves the Proposition.

\medskip

There is a similar result for 1-forms with values in (the additive 
group of) a vector space, - 
say, in the space of scalars $R$. The proof is  simpler, but similar. It 
is sketched in \cite{DK}, and it was one of the motivations for that paper.

\subsection{Primitives of closed group-valued 1-forms}Let $M$ be a manifold and 
$G=(G,*)$ a group. If $f:M\to G$ is a function, 
we get a $G$-valued 1-form\footnote{sometimes called 
the {\em Darboux derivative} of $f$; it is $f^{*}$ applied to the 
Maurer-Cartan form $a^{-1}b$ on $G$.} $df$ as follows: Let $x\sim y$ in $M$. Then we put
$$df(x,y):= f(x)^{-1}*f(y).$$
This is clearly a closed form. If $\omega$ is any $G$-valued 1-form 
on $M$, and $\omega =df$ for some $f:M\to G$, we say that $f$ is a 
{\em primitive} of $\omega$. So a necessary condition for $\omega$ to 
have a primitive is that $\omega$ is closed. If $U\subseteq M$ is 
``formally open'', (meaning: $x\in U$ and $x\sim y$ implies $y\in 
U$), then we may have a function $f:U \to G$ satsifying $df(x,y)= 
\omega (x,y)$ for $x\sim y$ in $U$, a primitive of $\omega$ on $U$. 
In global terms, it may be that 
$M$ is can be covered by such $U_{i}$s, and possessing  
primitives on each $U_{i}$ 
but with obstructions to 
patching these ``partial'' primitives together to a global function $f:M\to G$. 
However, on the 
formal level, we have the following construction. For $x\in M$, let 
$\M_{\infty}(x)$ be the set of points $y\in M$ for which there exists 
an $n$-path (for some $n$) with $x$ as domain and $y$ as codomain. 
This is clearly a formally open subset of $M$. 
Now the following is an easy Corollary of Proposition \ref{quadr+x}
\begin{cor}\label{primitivx}Let $\omega $ be a closed $G$-valued 1-form on $M$. Then 
for each $x_{0}$, there exists a unique partial primitive $f$ for 
$\omega$, defined on $\M_{ \infty} (x_{0})$ and with $f(x_{0})=1$.
\end{cor}
{\bf Proof.} 
Let $y\in \M_{ \infty} (x_{0})$, so there exists (for some $n$)  an 
$n$-path $\underline{x}$ from $x_{0}$ to $y$, say $(x_{0}, x_{1}, 
\ldots ,x_{n-1},y)$. We put $f(y):= \int _{\underline{x}}\omega$. By Proposition 
\ref{quadr+x}, this is, 
for given $n$, 
independent of the choice of the path. It may be that $y$ is the 
codomain of a shorter path, say of length $m<n$;  such path 
$\underline{z}$ may 
be augmented by a $n-m$ copies of $y$ in the codomain end, to provide 
an $n$-path $\underline{z}'=(\underline{z}, y, \ldots,y)$; but 
$\int_{\underline z} \omega =  \int_{\underline{z'}} \omega$, because 
$\omega (y,y)=1$. So $f$ is well defined. Furthermore $df = \omega$. 
For, if $y\sim z$  and if $y$ can be reached from $x=0$ by an 
$n$-path $\underline{x}$, then we may use $\underline{x}$ and 
the $n+1$-path $(\underline{x},z)$ to describe $f(y)$ and $f(z)$ 
respectively; and these two paths show that
$f(z) = f(y)*\omega (y,z)$, or equivalently, $\omega 
(y,z)=f(y)^{-1}f(z) (= df(y,z))$. So the constructed $f$ is indeed a primitive of 
$\omega$. 
The uniqueness of $f$ follows easily by induction in $n$.

\medskip
Note that the paths in $M$ form a category, by concatenation of 
paths; and that $\int \omega$ is takes composition in this category 
to multiplication $*$ in $G$.

\medskip

In the following, we assume that the manifold $M$ is path connected, 
meaning that any two points in $M$ can be connected by an $n$-path, 
for some $n$; equivalently, for all $x\in M$, we have 
$\M_{\infty}(x)=M$. This is a strong smallness condition; in fact, if 
$\sim$ is trivial in the sense that
$x\sim y$ implies $x=y$, then only one-point spaces are path connected! 
However, in well adapted toposes, and in algebraic geometry, the infinitesimal neighbour relation $\sim$ is 
not trivial. 

\medskip

\small

More importantly, though, is the fact that the 
infinitesimal constructions form a blueprint of what kind of 
approximations can be made in the physical world, where one considers 
small steps as  infinitesimal (building a round chimney out of 
square bricks, or forming Riemann sums). But it goes also, and 
primarily, the other way: from the small steps in the real world, one 
gets the geometric idea for rigorous $\sim$-infinitesimal notions - out of which 
even may grow rigourous analytic calculations (say in the form of 
power series). Often, the calculations is 
all you are presented with, as if they had  dropped from the skyes. 

\normalsize

\section {Connections in groupoids}\label{CGx}

The content of the present Section is presented in more detail in 
\cite{CPCG}.

We consider a groupoid $\Phi \rightrightarrows M$, where $M$ is 
equipped with a reflexive symmetric relation $\sim$. 
 Recall from \cite{E}, \cite{V} or \cite{SGM} that a connection in such groupoid 
may be defined as a map $\nabla : M_{(1)}\to \Phi$ with $\nabla 
(x,x)=1_{x}$ and $\nabla (y,x)= \nabla (x,y)^{-1}$.

The connection $\nabla$ is called {\em 
flat} (or {\em curvature-free}) if 
\begin{equation}\label{flax}\nabla (x,y).\nabla (y,z) = \nabla (x,z) ,\end{equation}
whenever $x\sim y$, $y\sim 
z$ and $x\sim z$, 
in analogy with (\ref{closedx}). (We compose from left to right in 
$\Phi$.)
In fact (\ref{closedx}) may be seen 
as the special case where the groupoid $\Phi \rightrightarrows M$ is 
$M\times M \times G$, and the connection is given by $\nabla (x,y):=  
(x,y, \omega (x,y))$ (a ``constant'' groupoid with vertex group $G$). For a groupoid which 
is locally of this form, one may locally choose such a trivialization,  and 
in terms of that, one can encode the connection by a $G$-valued 
1-form, which is closed iff the connection is flat.
 Therefore, for a flat connection $\nabla$ in such a groupoid, 
the Proposition \ref{QLx} implies that
\begin{equation}\label{flatgx}\nabla (x,y_{1}).\nabla (y_{1},z)= \nabla 
(x,y_{2}).\nabla (y_{2},z),\end{equation} for any  $\sim$ quadrangle 
$(x,y_{1}, y_{2},z)$ (meaning that $x\sim y_{i}\sim z$ for $i=1$ and 
$i=2$); and this statement does not depend of the choice of the local 
trivializations.

There is a more general notion of non-holonomous\footnote{More 
completely: ``not-necessarily-holonomous'' 
connections.} connection in such a groupoid; 
it is a law  which to an $n$-path 
$\underline{x}$ in $M$
 associates an arrow
$\nabla (\underline{x}): x_{0}\to x_{n}$ in 
$\Phi$; such laws are called {\em non-holonomous} connections (of 
order $n$), cf.\ 
considered in \cite{E} and \cite{V}; se also \cite{CNHJ}. Any 
connection $\nabla$ in $\Phi \rightrightarrows M$ gives rise to such 
non-holonomous connection of order $n$, 
namely: to $\underline{x}=(x_{0}, \ldots ,x_{n})$, one associates the composite arrow in 
$\Phi$,
$$\begin{diagram}x_{0}&\rTo ^{\nabla (x_{0},x_{1})}&x_{1}& \rTo 
^{\nabla(x_{1},x_{2})}&x_{2}& \; \;  \cdot \cdot \rTo^{\nabla 
(x_{n-1},x_{n})}  &x_{n}
\end{diagram}.$$
This non-holonomous connection is denoted $\nabla * \nabla * \ldots* 
\nabla$ ($n$ times), or $\nabla^{*n}$, cf.\ \cite{V}.  In case $\Phi \rightrightarrows M$ is the 
 groupoid
$M\times M \times G \rightrightarrows M$, $\nabla$ may be identified 
with a $G$-valued 1-form, and $\nabla^{*n}(\underline{x})$ may be 
identified with $\int_{\underline{x}}\omega$.

We assume, as in the beginning of Section 1, that $G$ admits an 
(auxiliary) multiplication preserving embedding into an algebra $(W,*)$ (short: ``$G$ is a matrix group'').

The following\footnote{I believe that it was first proved in \cite{V}, 
Theorem 7.} is now an immediate generalization of Proposition 
\ref{quadr+x}. 

\begin{prop}\label{holonomx}Assume that, locally, $\Phi \rightrightarrows M$ admits some 
isomorphisms (over $M$) with  groupoids of the form $M\times G \times M$ for some 
matrix group $G$ ; then  if $\nabla$ is flat, 
$\nabla^{*n}(\underline{x})$ only depends on $x_{0}$ and $x_{n}$.
\end{prop}
{\bf Proof.} The auxiliary isomorphism allows us to translate the data of $\nabla$ 
into a $G$-valued 1-form $\omega$, which is closed iff $\nabla $ is 
flat. Then Proposition \ref{quadr+x} shows the independence.

\medskip 

Note that such an auxiliary isomorphism of $\Phi$ with  $M\times G \times M$ is not 
intrinsic to the geometry; but since the conclusion of the 
Proposition does not mention this auxiliary isomorphism, the 
conclusion is intrinsic to $\nabla$ and $\Phi \rightrightarrows M$. 

\medskip

The law $\nabla ^{*n}$ satisfying the conclusion of the Proposition 
is then what \cite{E} and \cite{V} would call a {\em holonomous} 
$n$th order connection in the groupoid, meaning that its value on an 
$n$-path only depends on the endpoints of the path. 

\medskip

It is clear that (whether $\nabla$ is flat or not), the construction 
provides a functor from the category of paths in $M$ to the category 
(groupoid) $\Phi \rightrightarrows M$. Thinking of the category of paths as a 
formal version of the category of (Moore-) paths in $M$, this functor is in 
terminology from \cite{V} (see also 5.8 in \cite{SGM}), the {\em path 
connection} given by $\nabla$.

\section{Affine connections}
Affine connections, in the combinatorial sense of \cite{CTC}, 
may be seen (Subsection \ref{ACGX} below) as a particular case of groupoid 
valued connections, as discussed in the previous Section. 

An affine 
connection is a certain  
structure $\lambda$ on a set $M$, equipped with a symmetric reflexive relation 
$\sim$. Namely $\lambda$ is a partially defined ternary operation,
$(x,y,z)\mapsto [zxy]$ on $M$, which is defined whenever $x\sim y$ and 
$x\sim z$.\footnote{The notation $\lambda$ for the ternary operation 
was used in \cite{CTC} and \cite{SGM}; what presently is denoted $[zxy]$ was in 
loc.\ cit.\ denoted $\lambda (x,y,z)$; the discrepancy in the ordering 
of the arguments will not be relevant presently, since we here hardly 
ever 
supply the symbol $\lambda$ with arguments. Essentially, the $[zxy]$- 
notation goes back to \cite{P}.
}
The axioms are: a book-keeping axiom, and three equational 
axioms. The book-keeping axiom is that for all such $x,y,z$: 
\begin{equation}\label{bookx}[zxy]\sim y \mbox{ and } [zxy]\sim 
z,\end{equation}
which may be depicted by
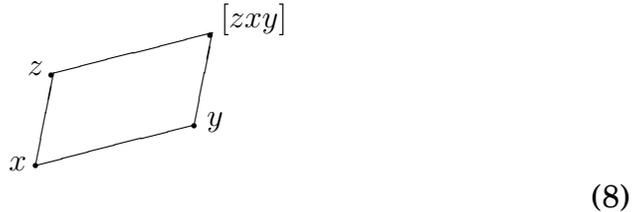
\begin{equation}
\begin{picture}(100,75)(0,-10)
\put(20,6){\line(4,1){60}}
\put(20,6){\line(1,5){7}}
\put(80,21){\line(1,5){7}}
\put(27,41){\line(4,1){60}}
\put(10,4){$x$}
\put(85,21){$y$}
\put(17,40){$z$}
\put(20,6){\circle*{2}}
\put(26,40){\circle*{2}}
\put(86,55){\circle*{2}}
\put(90,58){$[zxy]$}
\put(80,21){\circle*{2}}
\end{picture}
\label{aff-par}\end{equation}
in which the line segments display the $\sim$ relation.
The equational axioms are  two unit laws and one inversion law: the unit laws are
\begin{equation}\label{unitlx} [zxx]=z,\end{equation}
 \begin{equation}\label{unitrx} [xxy]=y.\end{equation}
 and the inversion law is
\begin{equation}\label{invrx} [[zxy]yx]=z.\end{equation}

The geometric meaning is that $[zxy]$ is 
the result of  
translating $z$ by that parallel translation which takes $x$ to $y$.
This process is in many models for the theory asymmetric in $y$ and 
$z$ ($z$ is ``passive'' (being moved), $y$ is ``active'' (is the mover); we shall here be interested in the case where further the 
 symmetry law holds:
\begin{equation}[zxy]=[yxz],\end{equation} in which case we call the affine connection {\em 
symmetric} (or {\em torsion free})
(and then the two unit laws of course are equivalent).

However, even without symmetry, the two unit laws and the inversion law suffice to construe 
an affine connection as a connection in the groupoid theoretic sense, 
as described in 
see Subsection \ref{ACGX} below.

\medskip
An affine connection $\lambda$ may be used for the following 
construction. Given two paths $\underline{y}$ and $\underline{z}$, with common domain $x$,  
say $\underline{y}=(x,y_{1},\ldots, y_{n})$ and 
$\underline{z}=(x,z_{1}, \ldots z_{m})$, we may form a  
2-dimensional $m\times n$ ``grid'' $u_{i,j}$ by induction:
We put $u_{0,0}:=x$ and $u_{0,j}:= y_{j}$, $u_{i,0}:=z_{i}$, and
\begin{equation}\label{2gridx}u_{i+1,j+1}:= 
[u_{i+1,j}u_{i,j}u_{i,j+1}].\end{equation}  
We call $x$ the {\em domain} of the grid, $u_{m,n}$ the {\em 
codomain}. We have, by the book-keeping laws for $\lambda$, that 
$u_{i+1,j}\sim u_{i,j}\sim u_{i,j+1}$. Note  that 
$u_{1,1}=[z_{1},x,y_{1}]$.
For this construction, we did not assume the symmetry law for 
$\lambda$. But if we also have symmetry of $\lambda$, it is clear 
that the construction of the 2-dimensional grid is likewise 
symmetric, in the sense that the grid, obtained by interchanging $\underline{y}$ and 
$\underline{z}$, is the transpose of the original grid (the $j,i$ 
entry in the transposed grid equals the $i,j$ entry in the original); 
so  in particular, the codomain of the 
grid ``spanned by'' $\underline{y}$ and $\underline{z}$ equals the 
codomain of the 
grid ``spanned by'' $\underline{z}$ and $\underline{y}$.

\subsection{Affine connections as groupoid connections}\label{ACGX}
Affine connections may be seen as a particular case of groupoid 
valued connections in the sense of Section \ref{CGx}.
Namely, for any manifold $M$, we have the groupoid $GL(M) 
\rightrightarrows M$, where an arrow $x\to y$ is a bijection $\M(x) 
\to \M (y)$ taking $x$ to $y$.\footnote{This groupoid is, for suitable notion 
of manifold, isomorphic to the locally constant groupoid consisting of fibrewise linear 
isomorphims $T_{x}(M) \to T_{y}(M)$, see Theorem 4.3.4 in \cite{SGM}, 
whence the choice of the acronym ``$GL$''. } 

More explicitly, for $x\sim y$ in $M$, the map $z\mapsto [zxy]$ 
defines a map $\widehat{\lambda}(x,y): \M(x) \to \M(y)$, by the (first) 
book-keeping law (\ref{bookx}), 
and it takes $x$ to $y$, by the unit law (\ref{unitrx}); it is a 
bijection with inverse $\widehat{\lambda}(y,x): \M(y) \to \M (x)$ by the inversion law 
(\ref{invrx}). It takes $x$ to $y$, by the other unit 
law (\ref{unitlx}). We shall also denote $\widehat{\lambda}$ by $\nabla$ (if $\lambda$ is 
understood), to conform with 
the notation of Section \ref{CGx},
$$\nabla (x,y) :=z\mapsto [zxy].$$

This viewpoint was likewise  introduced in \cite{CTC}, 
see also \cite{SGM} 2.3. The further requirement for $\lambda$, namely that $z\sim 
[zxy]$, 
we have here taken as a further book-keeping law (the second in 
(\ref{bookx})), even though it in 
synthetic differential geometry follows from general principles.

An affine connection $\lambda$ is called {\em 
flat}\footnote{Sometimes, e.g.\ in \cite{AM}, ``flat'' means what in 
our terminology is ``flat plus symmetric''.} if the corresponding 
groupoid valued connection $\widehat{\lambda}$ is flat (curvature free). Thus flatness 
 implies by the Quadrangle Law (Proposition \ref{QLx}) that 
transport of any $z\sim x$ around the two 2-paths in an arbitrary 
quadrangle with first vertex $x$ yield the same result. -  For the case where 
$\lambda$ is symmtric, we have some particular quadrangles which 
deserve the name  {\em \pgs}, namely quadrangles of the form $x,y,z, 
[zxy]$, as displayed in the picture (\ref{aff-par}); it
deserves the name: the {\em \pg}\; {\em spanned by} $y$ and $z$  with $x$
understood from the context; $x$ is called the 
{\em domain} or the {\em base} of the \pg; the {\em codomain} of the \pg \; is $[zxy]=[yxz]$).

We consider henceforth an affine connection $\lambda$  which is both symmetric 
and flat. 

The equation for moving $z\sim x_{0}$ (using $\widehat{\lambda}$) around the two 2-paths
from $x_{0}$ to $[x_{2}x_{0}x_{2}]$ in such \pg \;  gives same 
result, by 
flatness: 
\begin{equation}\label{weakflatx}[[zx_{0}x_{1}]x_{1}
[x_{1}x_{0}x_{2}]]=[[zx_{0}x_{2}]x_{2}[x_{2}x_{0}x_{1}]]\end{equation}
and is in simplified notation the equation (\ref{cub4x}) below.


The local 
triviality assumptions of Proposition \ref{holonomx} are   valid for 
the groupoid $GL(M)$, if $M$ is a manifold, (using charts from a 
vector space) and imply the following for a flat $\lambda$: for any path $x\sim 
y_{1}\sim y_{2}\sim \ldots \sim y_{n}$, and any $z\sim x$, 
the result $u_{n}$ of ``iterated transport of $z$ along the path''
\begin{equation}\label{zzx}[[[[z,x,y_{1}], y_{1},y_{2}], 
y_{2},y_{3}]\ldots ,y_{n-1},y_{n}]\end{equation}
is independent of the intermediate points $y_{1}, \ldots ,y_{n-1}$ (so
$\nabla ^{*n}$ is holnomous, in the terminology applicable for 
general groupoid valued connections).

\medskip



Suppose we are given two paths with domain $x$, say $\underline{y}$ 
and $\underline{z}$, as above,  and the resulting grid $\underline{y}\times 
\underline{z}$, with entries $u_{i,j}$ as constructed from it, as in 
(\ref{2gridx}). Using flatness of the affine
connection, we can then prove 
\begin{prop}\label{umnx}The point $u_{m,n}$ only depends on $x$, $y_{n}$, and 
$z_{m}$.
\end{prop}
{\bf Proof.} We have to prove that the point $u_{m,n}$ is independent 
of the choice of the paths $\underline{y}$ 
and $\underline{z}$. By symmetry, it suffices to prove that for fixed 
$\underline{z}$, it is independent of the choice of the path 
$\underline{y}$. This follows by induction in the length $m$ of the 
path $\underline{z}$. For $n=1$, this is a consequence of the 
flatness of $\lambda$: the $u_{n}$ in (\ref{zzx}) above is 
independent of the choice of $\underline{y}$, as we observed; and 
this $u_{n}$ is the one that in the grid $\underline{y}\times 
\underline{z}$ appears as $u_{1,n}$. The result now follows by 
applying the induction hypothesis to $z, z_{m}$ and $u_{1,n}$, 
using the bijection of paths from $x$ to $y_{n}$ on the one hand,  and paths from 
$z$ to $u_{i,n}$ (using transport along $xz$) on the other.

\medskip

Assume now that $M$ is path connected.
For given $x\in M$, we have therefore an everywhere defined binary 
operation $+_{x}$ given as follows:  $z+_{x}y$ is the codomain of the grid 
given by a path from $x$ to $z$ and a path from $x$ to $y$. This value 
does not depend on the paths chosen, by Proposition \ref{umnx}. Also, 
if $z\sim x \sim y$, we have
\begin{equation}\label{extensx}z+_{x}y=[zxy].\end{equation}
Since $M$ is assumed path connected, the transitive closure $\sim 
_{\infty}$ of $\sim$ is the trivial relation: for all $x$ and $y$ in 
$M$, we have 
$x\sim_{\infty}y$.
The following is then almost immediate:
\begin{prop}\label{32x}For fixed $x$, the binary operation $+_{x}$ is 
commutative, and has $x$ as a unit. The ternary operation $(x,y,z)\mapsto z+_{x}y$ defines
a symmetric affine connection, with respect to the trivial 
neighbour relation $\sim _{\infty}$ on $M$, and it  extends the given 
$\lambda$. 
\end{prop}
Note that we have not yet asserted associativity of $+_{x}$. This 
will be proved in Theorem \ref{36x} below.

The triviality of the relation $\sim_{\infty}$ means that we can forget about it, 
in particular, the book-keeping laws (\ref{bookx}) are trivially 
satisfied. We have by (\ref{extensx}) in fact extended the given 
affine connection 
$\lambda$, and may use the same notation  $[zxy]$ for this extended  
and   
everywhere defined operation. If we need to distinguish, we call the 
original connection the {\em small} (or $\sim$-restricted) one, the 
new extended we call the {\em big} (or unrestricted) one, and similarly 
for \pgs.

\subsection{The Cube Lemma}\label{CLX}
We come to the combinatorial core of this Section. We still consider 
a ($\sim$-restricted) affine connection $\lambda$  which is both symmetric and  flat. 
By 
symmetry of $\lambda$, we have a well defined notion of \pg, spanned 
by two neighbours $x_{1}$ and $x_{2}$ of $x_{0}$, and, more generally, we have a well 
defined 2-dimensional grid spanned by two paths with common domain 
$x_{0}$.

We now consider the case of {\em three} neighbours of $x_{0}$, and, 
more generally, of three 
paths with common domain $x_{0}$.

Given a point $x_{0}$, and three neighbour points $x,y,z$ of it.
Let us name these three  points $x_{1}$, $x_{2}$, and $x_{4}$, 
in some order.\footnote{The reason for choosing the name $x_{4}$, 
rather than $x_{3}$, will be given later.} We get three \pgs\; with base $x_{0}$ : 1) the one  
spanned $x_{1}$ and $x_{2}$, 2) the one  
spanned $x_{1}$ and $x_{4}$, and 3) 
the one 
spanned by $x_{2}$ and $x_{4}$.
These \pgs \; appear in the following picture as faces adjacent to 0 of the 
displayed cube (the point marked ``7'' will be argued after the 
calculation); for simplicity we have 
written $k$ for $x_{k}$ ($k=0,1,2,4$), and omitted commas.

Moving $4$ along the two paths from $0$ to $[102]=[201]$ give the same 
result, by the flatness of $\lambda$:
\begin{equation}\label{cub4x}[[401]1[102]]=[[402]2[201]];
\end{equation}
similarly moving 2, (or by renaming the three variables $x,y,z$, i.e.\ 
by permuting the indices 1,2,4)
\begin{equation}\label{cub2x}[[204]4[401]]=[[201]1[104]];
\end{equation}
and similarly, moving 1
\begin{equation}\label{cub1x}[[102]2[204]]=[[104]4[402]].
\end{equation}
The left hand side of (\ref{cub4x}) equals the right hand side of 
(\ref{cub2x}), by symmetry of $\lambda$; the
left hand side of (\ref{cub2x}) equals the right hand side of 
(\ref{cub1x}), by symmetry of $\lambda$; and the
left hand side of (\ref{cub1x}) equals the right hand side of 
(\ref{cub4x}), by symmetry of $\lambda$. Note that we have only been 
using flatness w.r.to \pgs\; (``weak flatness''). We conclude:
\begin{lemma}[Cube Lemma]Assume that $\lambda$ is a symmetric and 
(weakly) flat affine connection. Then all  six 
expressions appearing in the equations (\ref{cub4x}), (\ref{cub2x}) 
and (\ref{cub1x}) are equal.
\end{lemma}
This equal value is the point named ``
7'' in the following picture.

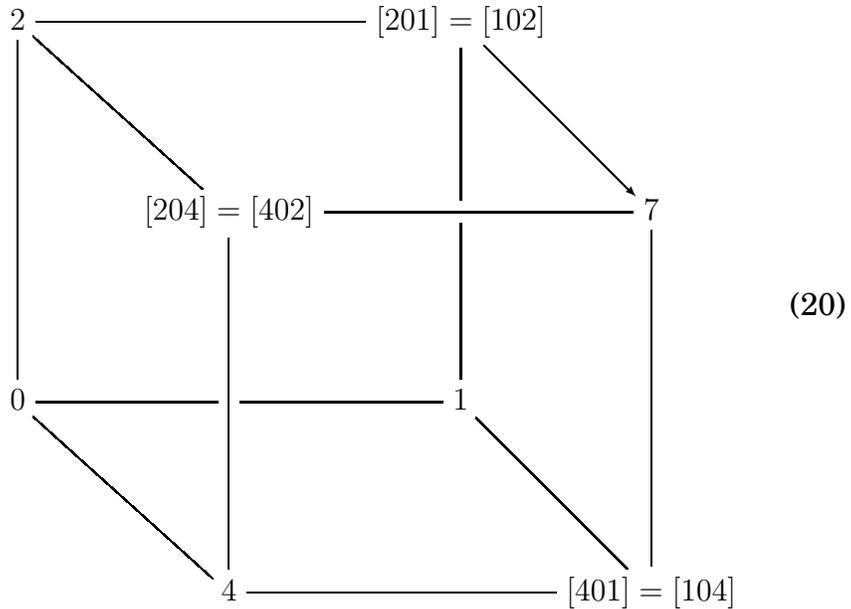
\begin{equation}\label{sevenx}\begin{diagram}
2&&\hLine&&[201]=[102]
&&\\
&\rdLine&&&\vLine&\rdTo & \\
\vLine&&[204]=[402]&\hLine&\HonV&&7\\
&&\vLine&&\vLine&&\\
0&\hLine& \VonH& \hLine&1&&\vLine\\
&\rdLine&&&&\rdLine&\\
&&4&&\hLine &&[401]=[104]\\
\end{diagram}\end{equation}

\medskip
 
\noindent (The naming of $x_{0}$ by $0$, and of $x$, $y$, and $z$ by $1,2,4$ (in some order) is a 
mnemotechnic device, with the purpose that the remaining points 
in the cube may be named $3,5,6,7$ in such a way that $[pqr]=p-q+r$; 
thus $[204]=6$, and $[623]=7$. 
The recipe for the naming is: consider the coordinate set of a vertex of the unit 
cube in ${\mathbb Z}^{3}$ as 
a number in digital notation; then write this number in decimal 
notation (just for compactness); e.g. the coordinate set of the point $[204]$ is 110 
which is digital notation for the number which in decimal notation is 
6. We invite the reader to write on the cube, writing  the ``3'' for 
$[102](=[201]$) etc.; ``7'' is then the 
equal value of any of the expressions in (\ref{cub4x}), (\ref{cub2x}) 
and (\ref{cub1x}).)

\medskip

(On the other hand, it is easy to see that if for a symmetric affine 
connection, the conclusion of the Cube Lemma holds, then this 
connection is weakly flat, i.e.\ the conclusion of the Quadrangle Law 
(Proposition \ref{QLx}) holds if the quadrangle is a \pg.)

\subsection{Three dimensional grid}
In the following Subsection, we shall strengthen the conclusion of 
Proposition \ref{32x} 
by adding a flatness assertion:
\begin{prop}If $\lambda$ is a symmetric flat  
($\sim$-restricted) affine connection $\lambda$  on a path connected $M$, then 
the  extension of $\lambda$ to an unrestricted affine connection is 
symmetric, and  flat with respect to (big) \pgs.
\end{prop}
{\bf Proof.} Only the flatness remains to be proved. The crux is to use the Cube Lemma for the restricted $\lambda$ to 
build a 3-dimensional grid (or ``big cube'')  $\underline{x}\times \underline{y}\times 
\underline{z}$   out of three paths  $\underline{x}$, 
$\underline{y}$,  $ \underline{z}$ with common domain, say 
$o$, and lengths $n$, $m$, and $k$, respectively. 
 The $i,j,l$ 
entry $w_{i,j,l}$ in the desired 3-dimensional grid is constructed  by 
induction, using the given affine connection $\lambda$. 
The initial conditionds are $w_{0,0,0}=o$, $w_{i,0,0}=x_{i}$, 
$w_{0,j,0}=y_{j}$, $w_{0,0,l}=z_{l}$.
The codomains of the three paths are denoted $x$, $y$, and $z$, 
respectively, thus $x=x_{n}$, $y=y_{m}$, $z=z_{k}$. 

The induction step uses crucially the Cube Lemma: 
$w_{i+1,j+1,l+1}$ is the last vertex  in the cube generated by 
$0:=w_{i,j,k}$, $1:= w_{i+1,j,k}$, $2:=w_{i,j+1,k}$, 
$4:=w_{i,j,k+1}$, as in the figure (\ref{sevenx}) (thus, $3= [ 
w_{i,j,k},w_{i+1,j,k},w_{i,j+1,k}]$ etc.)

The codomain  $w_{n,m,k}$ of the 3-dimensional grid is the point VII 
in the figure below.

Each of the six faces of the big cube is a 2-dimensional grid, and each of 
them can be seen as the witness of a \pg\;  for the unrestricted 
connection which we have constructed;   
thus the face containing $o,x,y$ is a grid constructing  $[xoy]$ for the 
unrestricted connection.
Similarly, the face containing the vertices $x$, $[xoy]$, $[xoz]$ is a 
grid constructing $[[xoy]x[xoz]]$; this is in the figure named 
``VII''. It also appears as a construction of other combinations 
of $o,x,y,z$, like $[[xoy]y[yoz]]$:

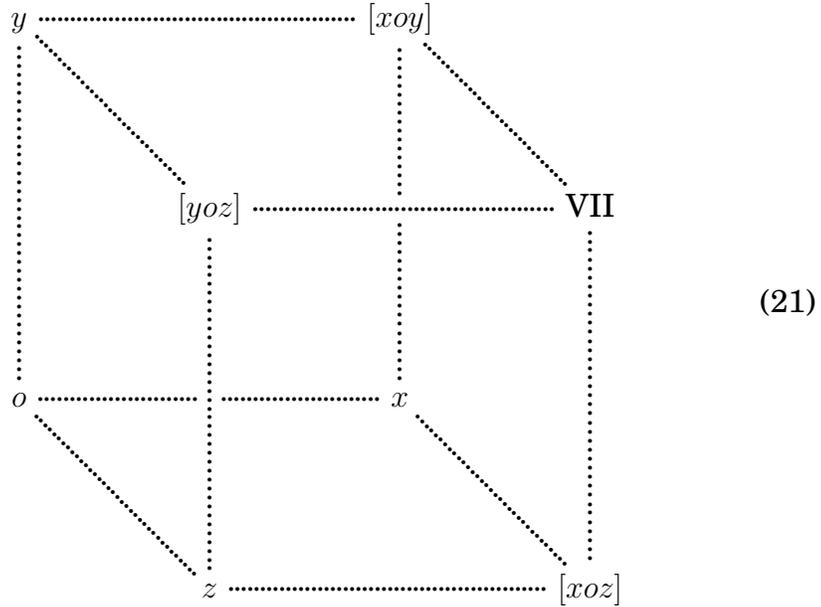
\begin{equation}\label{VIIx}\begin{diagram}
y &&\hDots&&[xoy]
&&\\
&\rdDots&&&\vDots&\rdDots & \\
\vDots&&[yoz]&\hDots&\HonV&&\mbox{VII}\\
&&\vDots&&\vDots&&\\
o&\hDots& \VonH& \hDots&x&&\vDots\\
&\rdDots&&&&\rdDots&\\
&&z&&\hDots &&[xoz]\\
\end{diagram}\end{equation}

\bigskip

\noindent The reader will observe that, except for the naming of the vertices, 
the cube in this figure looks like the cube in the previous one. But 
note the difference: in the previous one (the ``small'' cube), the lines indicate the 
$\sim$ relation, and the argument was that there were several 
constructions leading to the same result, which we then were allowed 
to give a name (choosing ``7'' for this name). 
In the present ``big'' cube, 
the lines indicate paths (therefore displayed as ``dotted'' lines) , where the 
three lines out of $o$ are 
arbitrary paths, and the rest of the cube is constructed canonically as 
the grid which these three paths generate.
The argument is now that the vertex VII (Roman notation for 7) is  {\em constructed} (as the last 
vertex $w_{n,m,k}$ of the grid), and we give 
{\em interpretations}
of it in terms of the unrestricted connection. More 
explicitly, the expressions in these equations (with suitable renaming)
express
the various ways we may see the, apriori existing, VII.

\medskip

There are some  interpretations of VII, available  in the big cube, 
whose analogs are not available for 7 in the small cube. These 
interpretations are based on the fact that one can concatenate paths. 
Thus, VII is the common value of the three expressions in (\ref{assx}):
\begin{prop}\label{cancx}[Cancellation law] For any $o,x,y,z$ on $M$, we have
\begin{equation}\label{cax}[[xoy]oz]  = 
[[xoy]y[yoz]]=[xo[yoz]].\end{equation} In particular, we have the 
associative law 
\begin{equation}\label{assx}[[xoy]oz]=[xo[yoz]].\end{equation}
\end{prop}
{\bf Proof.} The middle expression in (\ref{cax})  is the VII in the cube. To 
construct the right hand side, we must pick two paths:  from $o$ to $x$, 
and from $o$ to 
$[yoz]$. For the first, we pick the path $\underline{x}$, already 
used, and which appears in the grid as the path of points
$$o,w_{1,0,0}, w_{2,0,0},\ldots ,w_{n,0,0}.$$
For the second, we pick the concatenation of two of the paths that 
appear as edges in the cube; explicitly, the concatenated 
path is
$$o,w_{0,1,0}, \ldots ,w_{0,2,0}, \ldots ,w_{0,m,0}, w_{0,m,1}, 
\ldots , w_{0,m,k}.$$
It is clear that the 2-grid (of size $n\times (m+k)$) (constructed 
using the original (small) connection) has as its codomain $w_{n,m,k}$, which 
is also the codomain of the 3-grid used for VII.

\medskip

(Note that we cannot state the equation (\ref{cax}) for a general 
affine connection, since we do not have $[xoy]\sim o$, in general, so
the $\sim$-restriction  for forming
$[[xoy]oz]$ may not hold.)

\medskip 

Therefore, we have that
for fixed $o \in M$, the binary operation
$(x,y) \mapsto [xoy]$ is associative: both $[[xoy]oz]$ and $[xo[yoz]]$ 
equal $[[xoy]y[yoz]]$. If $M$ is path connected, this is an 
everywhere  defined binary operation $M \times M \to M$. It 
makes good sense  to denote this binary operation on $M$ by $x+_{o}y$:
$$x+_{o}y:=[xoy].$$
\begin{thm}\label{36x}For each $o\in M$ (assumed path connected), the binary operation $(x,y) \mapsto 
x+_{o}y$ makes $M$ into an abelian group. For any other $o'\in M$, 
the bijection $z\mapsto [zoo']$ is a group isomorphism.
\end{thm}
{\bf Proof.} We have already (Proposition \ref{32x}) that the operation $+_{o}$ is commutative, 
and also that $o$ is a unit. We just proved that it is associative. 
For existence of inverses, we have that $[oxo]$ will  serve: 
$$[xo[oxo]]= [[xoo]o[oxo]]=[[xoo]xo]=[xxo]=o$$
the second equality sign by the cancellation law, and the
first and the two last equality signs  by the unit law.

For the last assertion of the Theorem, we calculate:
$$[[xoo']o'[yoo']]=[[xoo']o'[o'oy]]=[[xoo']oy]$$
by symmetry and cancellation; and we continue:
$$= [[o'ox]oy]= [o'o[xoy]]=[[xoy]oo']$$
by  symmetry,  the associative law, and by symmetry again.
So if $f$ denotes the bijection considered, the total equation says 
that
$[f(x)o'f(y)]= f([xoy])$, and the Theorem is proved.

\medskip
The  Theorem here is really classical, going back to Pr\"{u}fer [P], 
who considered a ternary operation $\lambda$ satisfying similar equations as 
ours; but globally defined. He denoted  by $ 
(zx^{-1}y)$ what we denote $[zxy]$;  out of which he derives 
abelian group structures like $x+_{o}y$. 
In fact, his theory is, just as ours, an equational  presentation of the affine 
core of the theory of abelian groups, with this  ternary (but 
globally defined) 
operation as (the only) generator.

[P] calls a set with such a ternary operation a {\em Schar}; such 
structures, or generalizations thereof,  have been discovered 
indpendently by many authors, and under many names: by Baer, Certaine, 
Vagner, Lawson, and  others, including myself; 
see in particular Lawson's ``Generalised Heaps as Affine 
Structures'', in Hollings and Lawson
 {\em Wagner's Theory of Generalised Heaps}, 
Springer 2017, Cham), \cite{L}. 

In fact, a ``Schar'' or ``heap'' should just be termed: an {\em affine space over the 
ring ${\mathbb Z}$ of integers of scalars}.

What distinguishes our equations from Pr\"{u}fer's is that they admit 
restriction by a (reflexive symmetric) ``neighbour'' relation $\sim$. Some of the equations, 
valid in the unrestricted theory, do not make sense under such a 
restriction, see the remark after Proposition \ref{cax}. 

Theorem \ref{36x} may be reformulated, using known properties of affine 
spaces in general:


\begin{thm}\label{37x}Consider a manifold $M$ equipped with a symmetric and 
flat connection $\lambda$. Then every $\infty $-monad $\M$ in $M$ 
carries canonically structure of an affine space over ${\mathbb Z}$, with $\lambda 
(x,y,z)= y-x+z$, for any $y\sim x \sim z$ in $M$.
\end{thm}


\medskip

Our version is a ``formally  local'' one, i.e.\ for the notion of ``local'' derived from 
the notion ``formally open''. So in be coordinatized form, it gives 
only formal power series solutions, not anything about convergence. 
On the other hand, our result is canonical, 
whereas the classical result expresses that charts {\em 
exist} with certain properties, not the naturality of such charts.

\medskip

If we   replace the ring ${\mathbb Z}$ by the ring of reals ${\mathbb 
R}$, the (real) local result   is a version of a Theorem of Chern 1952, as 
quoted in \cite{AM} as   ``p.\ 108'' in these Chern notes.. 

\medskip

We shall sketch in the following Subsection how scalars (from 
${\mathbb R}$, say) may be introduced in the ${\mathbb Z}$-affine 
structure that we have constructed from $\lambda$.

\medskip

The conclusion of the (classical) result, as rendered synthetically in 
Theorem 3.7.4 in \cite{SGM} ( = Corollary 3 in \cite{CTC} = Theorem 2.3 in 
\cite{Buffalo}) imply this 
result, but note that in in these formulations, an abelian group (in 
fact, a vector space) is apriori {\em given} (the vector space on which the 
manifold $M$ is locally modelled), whereas in our formulation above, 
the abelian group is {\em constructed} (in the form of the affine structure 
on the $\infty$-monad). Furthermore, in loc.\ cit.\ we {\em assume} 
that a certain closed (group valued) 1-form is exact, which we here 
have {\em proved} in Corollary \ref{primitivx} to be the case. 

So what  the present paper adds to this classical 
theory is that it {\em derives} the global ternary operation, and the 
needed equations, out of infinitesimal data, namely the $\sim$-restricted 
(thus partially defined) affine connection.

\subsection{Affine combinations with scalars}\label{ACSx}

We take the notion of {\em affine space} over a commutative ring $R$ 
as meaning: ``affine combinations with 
coefficients from $R$, may be formed'' (recall that an affine combination is a 
linear combination where the sum of the coefficients is 
$1$). If we take $R$ to be ${\mathbb Z}$, we get the notion of 
``Schar'', ``heap'', ``commutative pregroup'' etc. To have more general coefficients (say 
${\mathbb Q}$), so that we e.g.\ can form the affine combination 
``midpoint'', $\tfrac{1}{2}x +  \tfrac{1}{2}y$, it suffices that we 
can form binary affine combinations, like 
$(1-t)x + ty$, for any $t\in R$, satisfying suitable compatibilities. 
We shall only be sketchy here.

\medskip

Let $R$ denote any commutative ring containing the rational numbers 
${\mathbb Q}$. We assume that we can form binary affine combinations 
with scalars from $R$, like $(1-t)x+ty$, for $x \sim y$. Such kind of 
structure we do have in the following two cases: 1) manifolds over 
$R$ (meaning: we can locally use charts from a KL vector space over 
$R$); and 2) general affine schemes $M$ over $R$, see \cite{ACAS}. The latter is purely 
formal, but the theory developed presently does not allow for this 
level of generality, since the groupoid $GL(M)$ is not necessarily locally 
constant, (say, if $M$ has singularities). For the manifold case, we 
have the technique of coordinate charts available, and hence the use 
of encoding the given affine connection in terms of Chrisoffel 
symbols:
$$\lambda (x,y,z)=[zxy]= z-x+y + \Gamma (x; z-x,y-x)$$
with $\Gamma$ bilinear in the arguments after the semicolon; and for 
$x\sim y$, the expression 
$(1-t)x+ty$ turns out not to depend on $\Gamma$ at all. See \cite{SGM}, 
2.3 for details. In fact, the monads $\M(x)$ and $\M(y)$ carry an 
action by the mutiplicative monoid of $R$, and the map 
$\widehat{\lambda}:\M(x) \to 
\M(y)$ induced by the affine connection $\lambda$ preserves the 
action, Proposition 2.3.7 in loc.\ cit.

 To indicate how the action of scalars extend to the whole of $M$ 
 (assumed path connected), I shall just indicate how to form $(1-t)x 
 + tu$ in case where $x$ and $u$ are ``second order'' neighbours, 
 i.e.\ in the case where there exists a 2-path $x\sim y \sim u$. Then 
 $u$ is of the form $[zxy]$ for some (unique( $z$ (take $z=[uyx]$). 
 Then we have (for $\lambda$ symmetric, equivalently, $\Gamma(x;-,-)$ 
 symmetric bilinear):
\begin{prop}  Let 
$u=\lambda (x,y,z)$.  
	Then $\lambda (x,y_{t},z_{t})$ only depends on $t$ and on $u$.
	\end{prop}
{\bf Proof.} In a coordinatized situation, let $y=x+d_{1}$ and 
$u=y+d_{2}$. Then $z= x+d_{2}- \Gamma (d_{1}, d_{2})$, where $\Gamma$ 
denotes the Christoffel symbol at the point $x$. Then $u= \lambda 
(x,y,z)$. The calculation for equation (3.3) in   
\cite{TCSx} 
gives that 
$$\lambda (x,y_{t},z_{t})= x+td_{1}+td_{2}-t\Gamma (d_{1},d_{2}) + 
\Gamma (td_{1},td_{2}),$$
and using that $\lambda$  is assumed symmetric, we have that $\Gamma 
$ is a symmetric bilinear for, so $\Gamma (v_{1},v_{2}) =\tfrac{1}{2}\Gamma 
(v_{1}+v_{2}, v_{1}+v_{2})$ for any pair of vectors $v_{1}$ and 
$v_{2}$ in $V$; thus with $y= x+d_{1}$, $z=x+d_{2}-\Gamma 
(d_{1},d_{2})$, 
$$\lambda (x,y_{t},z_{t})=x+t\cdot(d_{1}+d_{2}) + 
\frac{t^{2}-t}{2}\Gamma (d_{1}+d_{2}, d_{1}+d_{2}).$$
This clearly only depends on $t$ and $d_{1}+d_{2}$, i.e.\ on $t$ and 
$u= x+d_{1}+d_{2}$ as asserted.

\medskip
This means that we can define 
$(1-t)x + tu$ as $[x, (1-t)x+ty, (1-t)x+tz]$, independent of the 
``interpolating'' point $y$.

\medskip

Equational and foundational aspects of partially defined structures, like affine 
connection, with $\sim_{2}$ (like the above calculation) rather than 
$\sim _{1}$, may be found in \cite{FB}.

\medskip
Combining the (sketched) possibility of affine combinations with 
scalar coefficients with  Theorem \ref{37x}, we can state the 
following  

\begin{thm}Every flat and symmetric  affine connection locally comes about  from 
an actual affine structure, canonically constructed.
\end{thm}
(``locally'' in the sense of ``formally local'', i.e.\  on each $\infty$-monad).

\medskip

\small

\noindent The cubes were made with Paul Taylor's ``Diagrams'' package.

\medskip

\noindent Anders Kock

\noindent kock@math.au.dk

\noindent Dept.\ of Math., University of Aarhus,

\noindent Denmark

\noindent{February 2019}

\end{document}